\documentclass[pre,aps,10pt,superscriptaddress,notitlepage,nofootinbib]{revtex4-1}
\usepackage{epsfig,amsmath,amssymb,amsfonts,graphicx,color,calc,epstopdf}
\usepackage{stmaryrd}

 \let\t=\tau

\def\de{\mathrm d}

\def\CC{{\cal C}} 
\def\NN{{\cal N}} 
  
\def\DD{{\cal D}}

\def\ZZ{{\cal Z}}

\newcommand{\beq}{\begin{equation}} 
\newcommand{\eeq}{\end{equation}}
\newcommand{\ba}{\begin{eqnarray}}
\newcommand{\ea}{\end{eqnarray}}

\begin{document}

\title{Disordered high-dimensional optimal control}

  \author{Pierfrancesco Urbani}
  \affiliation{Universit\'e Paris-Saclay, CNRS, CEA, Institut de physique th\'eorique, 91191, Gif-sur-Yvette, France.}
 
 \begin{abstract}
Mean field optimal control problems are a class of optimization problems that arise from optimal control when applied to the many body setting. In the noisy case one has a set of controllable stochastic processes and a cost function that is a functional of their trajectories. The goal of the optimization is to minimize this cost over the control variables. Here we consider the case in which we have $N$ stochastic processes, or agents, with the associated control variables, which interact in a disordered way so that the resulting cost function is random. The goal is to find the average minimal cost for $N\to \infty$, when a typical realization of the quenched random interactions is considered. We introduce a simple model and show how to perform a dimensional reduction from the infinite dimensional case to a set of one dimensional stochastic partial differential equations of the Hamilton-Jacobi-Bellman and Fokker-Planck type. The statistical properties of the corresponding stochastic terms must be computed self-consistently, as we show explicitly. \end{abstract}
 
 \maketitle

\section*{Introduction}
Optimal control problems are ubiquitous in science and engineering. 
In a typical setting one has a dynamical system undergoing some deterministic or stochastic dynamics, over which the system has some degree of control. This is expressed through some set of control variables that the system can adjust along its dynamical trajectory. The goal of optimal control theory is to find the best control variables that minimize a designed cost function.

In a generic setting the solution of the problem is given in terms of a Bellman equation, which is a backward-in-time equation defined on a space whose dimensionality is inherited from the number of control variables that need to be fixed. When the dynamical evolution occurs in continuous time, the Bellman equation becomes a backward-in-time Partial Differential Equation (PDE) which is of the Hamilton-Jacobi type.
If the number of control variables is very large one has to deal with a PDE in high-dimension, whose solution is generically very difficult to find and to control. 
On the other hand, plenty of problems in statistical physics turn out to be exactly soluble in the high-dimensional limit; this relies on the fact that mean field theory becomes exact in high-dimension, provided that there is some degree of statistical equivalence between the degrees of freedom. When this is the case, the collective behavior of the full system is described by an effective problem in much lower dimension, whose properties must be fixed in a self-consistent way. In the following, we provide and example of this dimensional reduction applied to optimal control.

In recent years optimal control problems in high-dimension appeared in two interesting settings. 
The first one is the supervised learning problem for deep neural networks. In this case the goal is to fix the weights of the neural network in such a way that the latter performs a task, such as classification. 
If the network is deep (\emph{i.e.}, if the neurons are organized in a large collection of consecutive layers connected to each others), the classification process  can be thought of as a dynamical system in which the input is transformed progressively, layer by layer, into the output of the network. The dynamical evolution of the input signal can be tracked through the successive layers by a controlled dynamical equation, where the control variables are the weights of the network.
Therefore, the training dynamics can be thought as an optimal control problem. This point of view has been recently proposed in \cite{WJL19}. The corresponding optimal control problem is naturally high-dimensional in deep networks, in which the number of neural weights layer after layer is  large. Furthermore, the problem is disordered in nature since one wants to find the control variables that perform well the classification task when an \emph{ensemble} of inputs is considered.  As in plenty of other cases \cite{EV01}, one can think about the training set (\emph{i.e.}, the set of inputs used to train the network and fix its weights) as quenched disorder in the optimization problem.

Another setting in which high-dimensional optimal control problems appear naturally is that of Mean Field Games (MFG).
The theory of MFG is a rapidly expanding branch of analysis and partial differential equations \cite{GLL11, Ca10}. In a nutshell, one has a set of $N \gg 1$ agents whose state evolves in time. These agents are initialized in a given initial condition and aim at reaching a fixed goal at some final time, which may be infinite. Each agent has some control on their strategy to pursue the final goal, which determines its trajectory from the initial to the final time. Along their trajectory, the agents interact between themselves and are subjected to some potential cost due to the environment. Therefore, for each agent one has a cost function which depends on the collective position of the other agents during their time evolution. As discussed in the foundational works by Huang, Malham\'e, Caines \cite{HMC06} and Lasry and Lions  \cite{LL07}, if the interactions between agents are of the mean field type, one can perform a dimensional reduction in which the Nash equilibrium is described by a  solution of an optimal control problem for an effective agent which feels the effect of the other ones only through their density. This basic idea is of course very reminiscent of mean field theory in statistical physics. Mean field games with non-random interactions between the agents have been used to model several types of complex systems \cite{GLL11, GS14}, including active matter systems and in particular the flocking transition as described by the Cucker-Smale model \cite{CS07}. In several interesting applications the interactions between the agents are more complicated and heterogeneous, and thus we expect them to be suitably described by random terms.  While MFG with random environments have already been studied, see \cite{CDL16, CKZ20,De17}, the case of random interactions between the agents is largely open and, to the best of our knowledge, it is still to be investigated from the statistical physics perspective.

Motivated by these potential applications, in this work we study a high-dimensional, {prototypical} optimal control problem whose cost function contains disordered interactions between the control variables. We describe how to construct the associated mean field theory using techniques of disordered systems in statistical physics \cite{MPV87}. The main result of this analysis is a dimensional reduction in which one goes from a stochastic optimal control problem for a large number of interacting degrees of freedom, to an effective stochastic optimal control problem for a representative degree of freedom in a random environment whose statistical properties must be computed self-consistently. 
The paper is organized as follows. In Section \ref{sec:Model} we introduce the model.  
In Section \ref{sec:Conto} we set up the formalism, following closely \cite{Ka95,Ka05}. In Section \ref{sec:Replicas} we treat the model with the replica approach and show how the replica symmetric closure of the equations gives rise to an effective stochastic optimal control problem in a random environment, whose statistical properties are self-consistently determined. We conclude by discussing some perspectives on future directions and applications of these results.

\section{A simple model of disordered high dimensional optimal control}\label{sec:Model}
We consider a set of one dimensional\footnote{What follows can be easily generalized to higher dimensions.} random walkers obeying the controlled stochastic processes
\beq
\dot x_i= u_i(t) +\xi_i(t)
\eeq
where the noise is Gaussian and such that
\beq
\langle \xi_i(t)\rangle_\xi =0  \ \ \ \ \ \ \ \langle \xi_i(t)\xi_j(t')\rangle_\xi =\delta_{ij}\delta(t-t')
\label{noise}
\eeq
and we have denoted with brackets the average over the functional distribution of $\{\xi_i(t)\}$. 
The initial condition for these walkers may be deterministic, for example
\beq
x_i(0)=0
\label{init}
\eeq
but our computation can be extended to the case in which one has a separable probability distribution over the initial conditions, meaning that
\beq
P_0(\underline x(0)) = \prod_{i=1}^N \hat P_0(x_i(0))\:.
\eeq
In the following we consider the deterministic setting for simplicity and in particular Eq.~\eqref{init}, just to fix the ideas.
The variables $u_i(t)$ are time dependent and must be choosen in such a way that they minimize the following cost function
\beq
\CC[\underline x(0), 0] =\left\langle \frac{1}{2}\int_0^{t_f} \de \t \sum_{i=1}^Nu_i^2(\t) +\int_0^{t_f}\de \t V(\underline x(\t)) + \sum_{i=1}^N \phi(x_i(t_f)) \right\rangle_\xi
\label{cost}
\eeq
and we denoted with $\underline x(t) = \{x_1(t),\ldots, x_N(t)\}$.
The cost function in Eq.~\eqref{cost} is composed by three terms. 
The first one is an elastic running cost on the control variables $\underline u(t) = \{u_1(t), \ldots , u_N(t)\}$. Since the cost is quadratic in $u_i$, we call it quadratic mean field optimal control problem \cite{USG19}.
The second one contains the cost of interactions between agents and the last one is a separable terminal cost function that each agent tries to minimize. 
We consider the following interaction potential
\beq
V(\underline x) = \sum_{i=1}^{N} \nu(x_i) + \sum_{i<j} J_{ij} x_i x_j
\eeq
where the local potential $\nu(x)$ is confining enough that the random walkers are not allowed to escape to infinity.
A simple example could be
\beq
\nu(x) = \frac 12 x^2 + \frac 1{4!} x^4\:.
\eeq
We leave $\nu(x)$ unspecified in what follows.
The couplings between the walkers are randomly distributed and we choose them from a Gaussian distribution
\beq
J_{ij} = \frac{J}{\sqrt N} z_{ij} \ \ \ \ \ \ \ \ P(z_{ij})=\frac 1{\sqrt{2\pi}}e^{-z_{ij}^2/2}\:.
\eeq
We stress that the random couplings are such that we cannot write the cost function in terms of the density field of the agents. 
The terminal cost, the last term in Eq.~\eqref{cost}, has been assumed to be separable among agents.
An example could be
\beq
\phi(x) = \frac12 (x-1)^2
\eeq
so that we ideally want the agents to go from $x=0$ to $x=1$.
However our formalism will follow also for more complicated cases, and we leave the form of the terminal cost unspecified in the following.
The cost function of Eq.~\eqref{cost} is random because it depends on the $J_{ij}$s. 
Therefore we want to compute its average value over such couplings, meaning
\beq
r_0=\lim_{N\to \infty} \frac 1N \overline{\min_{\underline u}\CC[\underline x(0), 0] }
\eeq
and the overline stands over the average over the random variables $J_{ij}$s.
The solution to this problem can be constructed in the following way.
We define the optimal \emph{cost-to-go} function $f(\underline x,t)$ as the optimal cost for the walkers starting in $\underline x$ at time $t$ and walking up to time $t_f$ as
\beq
f(\underline x,t) = \min_{\underline u} \left\langle \frac{1}{2}\int_t^{t_f} \de \t \sum_{i=1}^Nu_i^2(\t) +\int_t^{t_f}\de \t V(\underline x(\t)) + \sum_{i=1}^N \phi(x_i(t_f)) \right\rangle.
\label{MB_HJB}
\eeq
Then $f(\underline x,t)$ satisfies the Hamilton-Jacobi-Bellman (HJB) equation \cite{BD15}
\beq
\begin{cases}
-\partial_t f(\underline x, t) = -\frac 12 |\nabla_x f(\underline x, t)|^2 +\frac 12 \nabla^2_xf(\underline x, t) + V(\underline x)\\
f(\underline x,t_f) = \sum_{i=1}^N\phi(x_i(t_f))\:.
\end{cases}
\eeq
The HJB equation must be solved backward in time and its boundary condition is provided by the final cost the agents want to minimize.
Once the solution is found we get the optimal cost as
\beq
r_0=\lim_{N\to \infty} \frac 1N \overline{f(\underline x(0),0)}\:.
\label{optimal_cost_AA}
\eeq
Furthermore one can show, see \cite{BD15}, that the optimal strategy is given by choosing 
\beq
\underline u(\underline x,t) = - \nabla_x f(\underline x,t)\:.
\label{opt_u}
\eeq
The difficulty comes from the fact that we are interested in the high-dimensional limit $N\to \infty$ and therefore we need to solve the HJB equation in that limit. However we will show that precisely because of the high-dimensional limit, the solution of the problem can be computed explicitly. In particular we will show how to perform a dimensional reduction that takes into account the average over disorder and  that allows us to compute the optimal cost in Eq.~\eqref{optimal_cost_AA}.

\section{The structure of the computation}\label{sec:Conto}
We closely follow \cite{Ka95, Ka05} and consider a logarithmic (Cole-Hopf) transformation \cite{Fl77} to define
\beq
f(\underline x,t) = -\log \psi(\underline x,t)
\eeq
so that we have
\beq
r_0=-\lim_{N\to \infty} \frac 1N \overline{\log \psi(\underline 0,0)}
\eeq
where we have assumed for simplicity the initial conditions as in Eq.~\eqref{init}.
Given the equation on $f$ one can show that $\psi(\underline x,t)$ satisfies the following equation
\beq
-\partial_t \psi(\underline x,t)  = \left[-V(\underline x) + \frac 12 \nabla^2_x\right] \psi(\underline x,t)\:.
\eeq
with boundary condition $\psi(\underline x,t)= e^{-\sum_i^N\phi(x_i)}$.
Defining $\rho(\underline y, t| \underline x,t')$ the solution of the following (forward in time) problem
\beq
\begin{cases}
{\partial_t}\rho(\underline y, t| \underline x,t') = \left[-V(\underline y) + \frac 12 \nabla^2_y\right]\rho(\underline y, t|\underline  x,t')\ \ \ \ \ \ \ t>t'\\
\rho(\underline x', t| \underline x,t) =\prod_{i=1}^N \delta(x_i' - x_i)
\end{cases}
\eeq
one can show \cite{Ka95} that 
\beq
\psi(\underline 0,0) = \int \de \underline y  \rho(\underline y, t_f|\underline 0,0)e^{-\sum_{i=1}^N\phi(y_i)}\:.
\eeq
It follows that there is a path integral representation for $\psi(\underline 0,0)$, see \cite{Ka95}, as 
\beq
\psi(\underline 0,0) = \int_{\underline x(0)=\underline 0} \DD \underline x \exp\left[-S[\underline x]   \right]
\eeq
where the action is given by
\beq
S[\underline x]  = \sum_{i=1}^N \phi(x_i(t_f)) +\frac 12 \sum_{i=1}^N \int_0^{t_f} \de \t\left[\dot x_i^2(\t)+ \nu(x_i(\t)) \right]+\sum_{i<j}J_{ij} \int_0^{t_f}\de \t x_i(\t)x_j(\t)\:.
\eeq
We note that this action is rather similar to the large $N$ random manifold in one dimension \cite{MP91}, with the main difference that the random potential, once averaged out, will couple different times as we will see below.
Therefore we arrive at the final equation
\beq
r_0=-\lim_{N\to \infty} \frac 1N \overline{\log \int_{\underline x(0)=\underline 0} \DD\underline  x \exp\left[-S[\underline x] \right]}\:.
\eeq
The path integral lives in high-dimension when $N\to \infty$ and therefore it is hard to compute numerically.
Furthermore we need to sample the $J_{ij}$s and compute the average over these random variables. We will address both problems using the replica method and exploiting the mean field structure of the model. 
The main difficulty is to compute the average over the disorder which involve the average of the logarithm, a notoriously hard problem. We use the replica method of statistical physics of disordered systems \cite{MPV87} to address it.

\section{The replica approach}\label{sec_replicas}\label{sec:Replicas}
We start by defining the partition function
\beq
Z=\int_{\underline x(0)=\underline 0} \DD \underline x \exp\left[-S[\underline x]   \right]\:.
\eeq
The replica method builds up on the identity \cite{MPV87}
\beq
\overline {\ln Z} =  \lim_{n\to 0} \partial_n \overline{Z^n}\:.
\eeq
If we assume that $n$ is integer then we can write the following expression for the powers of $Z$:
\beq
\begin{split}
 \overline{Z^n} &=  \overline{ \int
_{\{\underline x^{(a)}(0)=\underline 0\}_{a=1,\ldots ,n}}\left[\prod_{a=1}^n  \DD \underline x^{(a)}\right] \exp\left[-\sum_{a=1}^nS[\underline x^{(a)}]\right]   }\:.
\end{split}
\eeq
The main advantage of having integer $n$ is that we can perform explicitly the average over disorder. 
Once this is done one needs to find a way to analytically continue $n\to 0$.
Averaging over the random couplings $J_{ij}$ we get
\beq
\overline{\exp\left[\sum_{i<j}J_{ij} \int_0^{t_f}\de \t x_i(\t)x_j(\t)\right]} \simeq \exp\left[\frac{J^2}{4N}\int_0^{t_f}\de \t\int_0^{t_f}\de \t' \left(\sum_{i=1}^Nx_i^{(a)}(\t)x_i^{(b)}(\t')\right)^2\right]
\eeq
and we have neglected subextensive terms in the $N\to \infty$ limit.
Using the Hubbard-Stratonovich transformation we can write
\beq
r_0 = - \lim_{N\to \infty} \lim_{n\to 0} \frac {1}N  \partial_n \int\DD Q \exp\left[-\frac {NJ^2}4 \sum_{ab}\int_0^{t_f}\de \t \int_0^{t_f}\de \t' Q_{ab}^2(\t,\t') + N\ln \ZZ[Q]\right]
\label{r0_rep}
\eeq
where
\beq
\begin{split}
\ZZ[Q]&=\int_{\{ x_a(0)=0\}_{a=1,\ldots, n}} \left[\prod_{a=1}^n\DD x_a\right] \exp\left[-\sum_{a=1}^n \tilde S[x_a] + \frac{J^2}2 \int_0^{t_f}\de \t \int_0^{t_f}\de \t' \sum_{ab}Q_{ab}(\t,\t')x_a(\t)x_b(\t')   \right]\\
\tilde S[x_a] &= \phi(x_a(t_f)) + \int_0^{t_f} \de \t\left[\frac {1}2\dot x_a^2(\t)+ \nu(x_a(\t)) \right]\:.
\end{split}
\eeq
The measure $\DD Q$ in Eq.~\eqref{r0_rep} is properly normalized.
When $N\to \infty$ we can use the saddle point method to evaluate the functional integral.
The equation for the saddle point is  
\beq
Q_{ab}(\t,\t') = \langle x_a(\t) x_b(\t')\rangle_{\ZZ}
\label{rep_SP}
\eeq
where we have defined the average as
\beq
\langle O \rangle_\ZZ =\frac 1{\ZZ[Q]} \int_{\{ x_a(0)=0\}_{a=1,\ldots, n}} \left[\prod_{a=1}^n\DD x_a\right]\exp\left[-\tilde S[x_a] + \frac {J^2}2 \int_0^{t_f}\de \t \int_0^{t_f}\de \t' \sum_{ab}Q_{ab}(\t,\t')x_a(\t)x_b(\t')   \right] O
\eeq
At this point we need to find the solution of the saddle point equation \eqref{rep_SP} and find how such solution behaves when we take the analytical continuation down to $n\to 0$. 
We note that Eq.~\eqref{rep_SP} is very close to the one emerging from the replica treatment of the quantum partition function of mean field quantum spin glasses \cite{BM80}.
However in that case trajectories are periodic in (imaginary) time (meaning that $x(0)=x(t_f)$) and $t_f$ plays the role of the inverse temperature which is due to the fact that one is trying to compute the trace of the Boltzmann factor using the Suzuki-Trotter formula. In that case one can show that translational invariance on the torus gives rise to the fact that $Q_{a\neq b}(\tau,\tau')$ does not depend on $\t$ and $\t'$. 
Here the situation is different since we do not have any periodic condition on trajectories. 
Therefore on very general grounds we only have
\beq
Q_{aa}(\tau,\tau')=Q_{aa}(\tau',\tau)\ \ \ \ \ \ \ Q_{ab}(\tau,\tau')=Q_{ba}(\tau',\tau)\:.
\eeq
To move forward, we need to consider an ansatz for the form of the solution $Q_{ab}(\t,\t')$ which is amenable to the $n\to 0$ analytic continuation.

\subsection{The replica symmetric ansatz}
We assume that the solution of the saddle point equations is of the replica symmetric form
\beq
Q_{ab}(\t,\t') = D(\t,\t') \delta_{ab} + (1-\delta_{ab})F(\t,\t')\:.
\eeq
With this assumption one can rewrite $\ln \ZZ$ in terms of functional Gaussian integrals as
\beq
\begin{split}
\ln \ZZ &= \ln \int \DD H \exp\left[-\frac 12 \int_0^{t_f}\de \t \int_0^{t_f}\de \t' H(\t)F^{-1}(\t,\t')H(\t')\right] \hat \ZZ^n(H)\\
&\simeq n \int \DD H \exp\left[-\frac 12 \int_0^{t_f}\de \t \int_0^{t_f}\de \t' H(\t)F^{-1}(\t,\t')H(\t')\right] \ln \hat \ZZ(H)
\end{split}
\eeq
where we assumed that $n\to 0$ and we have retained the leading term in $n$. Furthermore we used that the measure $\DD H$ is properly normalized, so that
\beq
 \int \DD H \exp\left[-\frac 12 \int_0^{t_f}\de \t \int_0^{t_f}\de \t' H(\t)F^{-1}(\t,\t')H(\t')\right] =1\:.
\label{norm_ga}
\eeq
Finally $F^{-1}$ is the inverse operator of the kernel $F$ and 
\beq
\hat \ZZ(H) =\int \DD h \exp\left[-\frac 12 \int_0^{t_f}\de \t \int_0^{t_f}\de \t' h(\t)(D-F)^{-1}(\t,\t')h(\t')\right] \hat \psi[h,H]\:.
\eeq 
Again the measure $\DD h$ is properly normalized analogously to Eq.~\eqref{norm_ga} and we have defined
\beq
\hat \psi[h,H] = \int_{x(0)=0} \DD x \exp\left[- \phi(x_i(t_f)) -  \int_0^{t_f} \de \t\left[\frac 12\dot x^2(\t)+ \nu(x(\t)) +J(h(\t)+H(\t))x(\t) \right] \right]\:.
\eeq
Therefore one can consider both $h(t)$ and $H(t)$ as uncorrelated Gaussian random fields with zero average and two point functions given by
\beq
\langle H(t)\rangle_H=0\ \ \ \ \langle h(t)\rangle_h=0 \ \ \ \ \langle H(t)H(t')\rangle_H=F(t,t')\ \ \ \ \langle h(t)h(t')\rangle_h=D(t,t')-F(t,t') \:.  
\eeq
We can then rewrite the saddle point equations in the following form
\beq
\begin{split}
D(\t,\t')&= \left\langle\frac{\langle \frac 1{J^2}\frac{\delta^2}{\delta h(\t) \delta h(\t')}\hat\psi[h,H]\rangle_h}{\langle \hat \psi[h,H]\rangle_h} \right\rangle_H\\
F(\t,\t') &=   \langle R(\t)R(\t') \rangle_H\\
R(\t)&=\frac{\langle \frac{1}{J}\frac{\delta}{\delta h(\t)}\hat\psi[h,H]\rangle_h }{\langle \hat\psi[h,H]\rangle_h}\:.
\end{split}
\label{SP_RS}
\eeq
We will now give an alternative route to compute both $\hat \psi[h,H]$ and the two point correlation functions that define the kernel $D$ and $F$, without resorting to the calculation of the full path integral.
It is easy to show, just by rewinding the same procedure we used to get to $\psi(\underline x,t)$, that one can write
\beq
\hat \psi[h,H] = e^{-c(0,0|h,H)}
\eeq 
where $c(x,t|h,H)$ satisfies the backward stochastic HJB equation
\beq
\begin{cases}
-\partial_t c(x,t|h,H) = -\frac 12 (\partial_x c(x,t|h,H))^2 +\frac 12 \partial_x^2c(x,t|h,H) + \nu(x) + J(h(t)+H(t))x\\
c(x,t_f|h,H) = \phi(x(t_f))\:.
\label{stoc_bel}
\end{cases}
\eeq
This equation is associated to a single body random stochastic optimal control problem. We have an effective one dimensional walker $\dot x(t) = u(t) + \xi(t)$, where $\xi(t)$ is a noise with the same properties in Eq.~\eqref{noise}, that tries to minimize
\beq
\CC_{h,H}(x',t)= \left\langle \frac{1}{2}\int_t^{t_f} \de \t u^2(\t) +\int_t^{t_f}\de \t \left[\nu(x(\t)) +J(h(\t)+H(\t))x(\t) \right] + \phi(x(t_f)) \right\rangle_\xi
\eeq
and we need to impose that the stochastic process starts at $x(t)=x'$.
The optimal strategy to solve the optimal control problem reads
\beq
u(x,t) = -\partial_xc(x,t|h,H)
\label{stoc_opt_prot}
\eeq
which is analogous to Eq.~\eqref{opt_u}.
We note that by taking Eq.~\eqref{stoc_bel} and deriving it with respect to $x$ one gets an equation for $u$ that is a forced stochastic Burgers equation with two stochastic terms given by the random fields $h$ and $H$ whose statistical properties are self consistently determined.
When the stochastic process for the effective agent is run with the optimal strategy provided by Eq.~\eqref{stoc_opt_prot} one can compute the associated Fokker-Planck equation
\beq
\begin{cases}
\partial_t \pi_{h,H}(x,t|y,t') = \frac 12 \left[ \partial_x^2\pi_{h,H}(x,t|y,t') + 2 \partial_x \left(\pi_{h,H}(x,t|y,t')\partial_xc(x,t|h,H)\right)\right]\\
\pi_{h,H}(x,t'|y,t')=\delta(x-y)
\end{cases}
\label{stoc_FP}
\eeq
This equation must be solved forward in time.
We note that $\pi_{h,H}(x,t|y,t')$ is a random measure since it depends on the full stories of $h(t)$ and $H(t)$ through $c(x,t|h,H)$.
Using again the result of \cite{Ka95} one can write
\beq
\hat\psi[h,H] = \int \de y \rho_{h,H}(y,t_f|0,0)e^{-\phi(y)}
\eeq
and $\rho_{h,H}(y,t_f|0,0)$ satisfies the following equation
\beq
\begin{cases}
{\partial_t}\rho_{h,H}( y, t| x,t') = \left[-\nu(y) -J(h(t)+H(t))y + \frac 12 \partial_y^2\right]\rho_{h,H}( y, t| x,t')\\
\rho_{h,H}( x', t| x,t) = \delta(x' - x)
\end{cases}
\label{Eq_rho_single}
\eeq
In order to compute numerically $\rho_{h,H}( y, t| x,t')$ on a given realization of the random fields $h$ and $H$, one can follow directly \cite{Ka95,Ka05} and consider a population of independent one dimensional walkers starting all at $x$ at time $t$ and all subjected to the same realization of $h(t)$ and $H(t)$.  They undergo pure diffusion and are killed (namely dropped from the simulation) with a rate that is $\nu(x)+J(H(t)+h(t))x$.
Here we also give an alternative way to compute $\rho_{h,H}( y, t| x,t')$. We can introduce $\psi_{h,H}(x,t)$ that satisfies
\beq
\begin{cases}
-\partial_t \psi_{h,H}(x,t) = \left[-\nu(x) -J(h(t)+H(t))x + \frac 12 \partial_x^2\right]\psi_{h,H}(x,t)\\
\psi_{h,H}(x,t_f) = e^{-\phi(x)}
\end{cases}
\label{eq_psi_single}
\eeq
and we have $\hat \psi[h,H] = \psi_{h,H}(0,0)$.
Using Eq.~\eqref{Eq_rho_single} and Eq.~\eqref{eq_psi_single} we can show that
\beq
\rho_{h,H}( x, t| y,t') =  \pi_{h,H}(x,t|y,t') \frac{\psi_{h,H}(y,t')}{\psi_{h,H}(x,t)}\:.
\eeq
Furthermore we define
\beq
\begin{split}
\llbracket O[h]\rrbracket &= \frac 1{\NN_h}  \left \langle O[h]e^{-c(0,0|h,H)}\right \rangle_h\\
\NN_h&=\left \langle e^{-c(0,0|h,H)}\right \rangle_h\:.
\end{split}
\eeq
With these definitions the replica symmetric equations are given by
\beq
\begin{split}
D(\t,\t') &= \left \langle  {\left \llbracket C(\t,\t')\right \rrbracket}    \right\rangle_H\\
F(\t,\t') &= \left \langle  {\left \llbracket m(\t)\right \rrbracket}{\left \llbracket m(\t')\right \rrbracket}  \right\rangle_H
\end{split}
\label{slef_cons}
\eeq
where we have defined 
\beq
\begin{split}
C(\t,\t') &= \int_{-\infty}^\infty \de x \int_{-\infty}^\infty \de x'\, xx'\pi_{h,H}(x,\t |x',\t')\pi_{h,H}(x',\t' |0,0)\\
m(\t)&=\int_{-\infty}^\infty \de x\, x\pi_{h,H}(x,\t |0,0)
\end{split}
\label{final}
\eeq
and we have assumed $\t\geq \t'$. 
This concludes the replica symmetric solution of the model.
The interpretation of these equations is rather clear.
We have reduced the problem of computing the solution of a deterministic and disordered HJB equation in infinite dimension, to the solution of a one dimensional stochastic HJB equation. In other words, thanks to the mean field nature of the interactions between the agents one has an effective agent that feels the interaction with the others through the random fields $h(t)$ and $H(t)$. The statistical properties of these random fields must be computed self-consistently through Eq.~\eqref{final}.
These equations tell that these random fields have both Gaussian statistics and one needs only to control their two point function.

Finally we can write the optimal cost as 
\beq
r_0= \frac{J^2}{4}\int_{0}^{t_f}\de \t\int_{0}^{t_f}\de \t'\left[D^2(\t,\t')-F^2(\t,\t')\right]-\langle \ln \NN_h \rangle_H
\label{final_cost}
\eeq
which concludes the replica symmetric treatment of the model.

We finally write an ideal algorithmic strategy to solve the saddle point equations.
One starts from a guess for the kernels $D(\t,\t')$ and $F(\t,\t')$.
Given these two kernels one extracts samples of the trajectories of $h(t)$ and $H(t)$ from their Gaussian probability distributions.
For each of these samples one computes the solution of the stochastic backward HJB equation, Eq.~\eqref{stoc_bel} and given a solution of this equation one gets the solution for the associated Fokker-Plank equation, Eq.~\eqref{stoc_FP}. The update of the kernels $D(\t,\t')$ and $F(\t,\t')$ is done using Eqs.~\eqref{slef_cons} and \eqref{final}.

We note that when the interactions between the agents vanish, $J\to 0$, one gets the usual optimal control problem for independent agents. $C$ and $m$ are then related to the dynamical two point functions (connected and disconnected) for a typical agent. However, since $J\to 0$, the final cost in Eq.~\eqref{final_cost} does not depend on these two point functions as it should.

\section{Conclusions and Perspectives}
We have introduced a simple model of a high-dimensional disordered optimal control problem where the running cost depends on a set of random couplings describing heterogeneous interactions between the agents.
We have presented the analysis of the model using the replica method under the replica symmetric ansatz, and we have derived a set of stochastic HJB and Fokker-Planck equations describing the optimal 
strategy of an effective agent. 
The stochastic equations must be solved self-consistently, {similarly to what happens commonly when studying the} dynamical mean field theory of disordered systems \cite{EO92,GKKR96,MKUZ20}. 

Let us comment on the connection between our results and the applications of high-dimensional optimal control mentioned in the introduction. A statistical physics approach to game theory has already been developed, even in presence of disorder. A prominent example of this are minority games \cite{MC01, CMZ05, Sh06}. In these cases one reaches a description of the problem in terms of an effective agent, in analogy with what we obtain here. The approach we adopt is however quite different, and more rooted into optimal control theory. It is natural to expect that the optimal control problem we discussed can be rephrased as a MFG once one considers the cost function of each agent \emph{given} the others playing optimally. One would then get a system of coupled infinite dimensional PDEs and one could conjecture that the self-consistent HJB equation we derived in this work describes the statistical properties of a typical agent in the game. We believe that the cavity method of disordered systems would be instrumental to prove explicitly this connection: we leave this to future work. For what concerns the supervised learning problem, we stress that our approach is also quite different with respect to the one of \cite{WJL19}, since in \cite{WJL19} the disorder is in the training set (which corresponds to a set of input-output associations of the dynamics) and the control is performed over an entire set of inputs. In the problem that we discussed, the disorder is instead in the cost function and does not depend on the initial condition of the dynamics. As a consequence, the resulting mean field theory is rather different.

The strategy of the computation we presented can be {readily} extended to more complex cases. For example, one can consider a situation in which the interaction between the agents
is multi-body, {generalizing with respect to the two-body interactions discussed above}.  
An interesting additional perspective of this work is the analysis of the solution of the stochastic PDEs that we have obtained when 
the cost function contains non-analyticities. In this case one could think that the kernels $D$ and $F$ have a scaling behavior to regularize the solution of the HJB close to the non-analytic points.
An example of this kind could arise when the phase space for the trajectories of the stochastic processes contains obstacles
and inaccessible regions which could create bottlenecks resulting in crowding effects. 
For instance, one may think to generalize the problem to agents moving in higher dimensions and to place forbidden regions along their trajectories.
Close to the boundary of these regions one could have interesting behaviors. 

Finally we did not attempt a numerical solution of the saddle point equations, which is left for future work. This goes together with the analysis of the stability of the replica symmetric assumption that we used to derive the effective single agent optimal control problem, and whose validity must be checked self- consistently.

\section{Aknowledgments}
This work was supported by ``Investissements d'Avenir" LabEx PALM (ANR-10-LABX-0039-PALM).
The author warmly thanks Cesare Nardini, Denis Ullmo and Valentina Ros for very useful discussions.

%


 \end{document}